\documentclass[nonblindrev]{informs3}

\usepackage{endnotes}
\usepackage{mathtools}
\let\footnote=\endnote
\let\enotesize=\normalsize
\def\notesname{Endnotes}%
\def\makeenmark{$^{\theenmark}$}
\def\enoteformat{\rightskip0pt\leftskip0pt\parindent=1.75em
  \leavevmode\llap{\theenmark.\enskip}}

\TheoremsNumberedThrough     
\ECRepeatTheorems

\EquationsNumberedThrough    

\usepackage{bm}
\usepackage{amsmath}
\usepackage{mathtools}
\newcommand{\ud}{\,\mathrm{d}}
\usepackage{algorithm}
\usepackage{algorithmic}
\usepackage{subcaption}
\usepackage{mathtools}

\usepackage{hhline}
\usepackage{booktabs}
\usepackage{bbm}
\usepackage{courier}
\usepackage{adjustbox}
\newcommand\numberthis{\addtocounter{equation}{1}\tag{\theequation}}

\begin{document}


\RUNAUTHOR{}

\RUNTITLE{}

\TITLE{Discrete Gradient Flow Approximations \\ of High Dimensional \\ Evolution Partial Differential Equations \\ via Deep Neural Networks}

\ARTICLEAUTHORS{%
\AUTHOR{Emmanuil H.~Georgoulis}
\AFF{1. School of Computing and Mathematical Sciences,
	University of Leicester, United Kingdom, \\ 
	2. Department of Mathematics, School of Applied Mathematical and Physical Sciences, \\ National Technical University of Athens, Greece, \\
	3. IACM-FORTH, Greece \\ \EMAIL{georgoulis@math.ntua.gr}  \URL{}}
\AUTHOR{Michail Loulakis}
\AFF{1. Department of Mathematics, School of Applied Mathematical \& Physical Sciences, \\ National Technical University of Athens, Greece,  \\ 2. IACM-FORTH, Greece\\ \EMAIL{loulakis@math.ntua.gr} \URL{}}
\AUTHOR{Asterios Tsiourvas}
\AFF{Operations Research Center, Massachusetts Institute of Technology, Cambridge, MA, USA,\\ \EMAIL{atsiour@mit.edu} \URL{}}
} 

\ABSTRACT{%
	We consider the approximation of initial/boundary value problems involving, possibly high-dimensional, dissipative evolution partial differential equations (PDEs) using a deep neural network framework. More specifically, we first propose discrete gradient flow approximations based on non-standard Dirichlet energies for problems involving essential boundary conditions posed on bounded spatial domains. The imposition of the boundary conditions is realized weakly via non-standard functionals; the latter classically arise in the construction of Galerkin-type numerical methods and are often referred to as ``Nitsche-type'' methods. Moreover, inspired by the seminal work of Jordan, Kinderleher, and Otto (JKO) \cite{jko}, we consider the second class of discrete gradient flows for special classes of dissipative evolution PDE problems with non-essential boundary conditions. These JKO-type gradient flows are solved via deep neural network approximations. A key, distinct aspect of the proposed methods is that the discretization is constructed via a sequence of residual-type deep neural networks (DNN) corresponding to implicit time-stepping. As a result, a DNN represents the PDE problem solution at each time node. This approach offers several advantages in the training of each DNN. We present a series of numerical experiments which showcase the good performance of Dirichlet-type energy approximations for lower space dimensions and the excellent performance of the JKO-type energies for higher spatial dimensions.
}%


\maketitle

\section{Introduction}

\par The numerical approximation of high-dimensional PDEs via collocation-type approaches using artificial neural network (ANN) architectures has received considerable interest in recent years \cite{sirignano_dgm_2018, PINN, e_deep_2017, jentzen, liao_deep_2019}. This activity appears to have been, at least partially, fuelled by the breakthroughs in the use of DNNs in multiple research domains \cite{goodfellow, lecun2015deeplearning}. Nonetheless, in the context of scientific computing at large, the use of artificial neural networks is yet in its early stages. 

\par In the context of numerical solution of partial differential equations (PDEs), multiple neural networks and deep learning approaches have been proposed in recent years. We discuss some recent developments that are close in spirit to the contribution below. In \cite{e_deep_2017} a deep learning minimization of the Dirichlet energy is proposed for the numerical approximation of elliptic PDEs, while in \cite{sirignano_dgm_2018}, the, so-called \emph{deep Galerkin method} (DGM) (and the related \emph{physics informed neural network} (PINN) approach \cite{PINN}) advocate the minimization of $L_2$-norm residuals of the PDE and the boundary/initial conditions in a collocation setting.  Further, in  \cite{liao_deep_2019}, the authors consider an alternative loss functional implementing weak imposition of boundary conditions in the classical Nitsche fashion from finite element analysis \cite{nitsche_uber_1971}, to alleviate the difficulties observed in \cite{e_deep_2017} to enforce conformity in the solution space. In \cite{wang_implicit_2020}, the difference between the Ritz-Galerkin method and the deep learning methods are studied with a focus on the implicit regularity that the second method imposes. In \cite{sheng_pfnn_2020} a penalty-free neural network approach for solving a class of second-order boundary value problems on complex geometries is presented that uses a combination of two neural networks, one for the domain and one for the boundary, to solve the problem. Also, in \cite{grohs_deep_2020} the authors propose a deep neural approximation for high-dimensional elliptic PDEs with boundary conditions for finite domains.


This work is primarily concerned with the development of a class of DNN-based collocation methods for dissipative evolution PDEs, as well as some performance comparison with known alternatives. The basic idea proposed in this work is the use of loss functions stemming from discretized gradient flows \cite{savare}. The proposed methods operate in a time-stepping fashion, by seeking to minimize the DNN parameters for each time interval. To showcase the potential generality of this approach, we consider heat-type equations viewed as  $L_2$-gradient flows of the Dirichlet energy and $2$-Wasserstein metric gradient flows of the Gibbs-Boltzmann entropy in the spirit of the seminal work of Jordan, Kinderleher and Otto (JKO) \cite{jko}. A series of numerical experiments using a variant of the ResNET-type architecture from \cite{sirignano_dgm_2018} highlight the good performance of the proposed methods. A comparison of the new scheme with the space-time residual approach of the Deep Galerkin Method from \cite{sirignano_dgm_2018} is also presented. The use of the JKO functional was advocated in \cite{asterios},  from which the present work originates.  We also note the recent manuscript \cite{new_JKO} proposing a fully-connected DNN architecture for various gradient flow functionals.

The new ANN-based collocation time-stepping methods proposed in this work use a different ResNET (with the same architecture) for representing the PDE approximation at each time step. This offers several theoretical and practical advantages compared to `monolithic' space-time residual minimization methods \cite{sirignano_dgm_2018, jentzen, PINN}, at least for low-to-medium dimensional evolution problems. More specifically, the ``time-instance" ResNETs can be chosen to have a significantly smaller width and depth, compared to full space-time collocation approaches, resulting in both faster performance and better accuracy: this is supported by numerical comparisons presented below. Moreover, the smoothness of the evolution semigroup facilitates the use of the optimized ResNET parameters at time-step $n$ to be used as an excellent initial guess for the optimization taking place at the next timestep $n+1$, thereby reducing significantly the number of epochs required in total. This is a particularly salient feature as the optimization steps initiate with the approximation of the known initial condition of the PDE problem (timestep $0$), i.e., a simple function fitting step not involving differentiation of the ResNET; this is typically achieved to very high accuracy. Having a very accurate representation of the initial condition gives rise to excellent starting values for the optimization to find the approximation at timestep $1$ and so on.  The superior performance of the proposed method is showcased through a series of comparative numerical experiments with known methods. 

To arrive at the above complete numerical framework, several developments are also considered for the various individual terms in the loss functions. In particular, for the $L_2$-norm gradient flow method, we consider weakly imposed Dirichlet-type boundary conditions through a consistent minimization functional in the spirit of Nitsche  (cf.~also \cite{liao_deep_2019}, as well as \cite{makridakis} for the same functional in the context of interface modeling in structural mechanics). A crucial methodological development proposed in this work, which may be of independent interest, is the method of computation of the unknown \emph{penalty parameter} required in the weak imposition of Dirichlet boundary conditions. More specifically, in the context of DNNs, the unknown penalty parameter is impossible to estimate directly via so-called inverse estimates as  is done in classical finite element methods. To alleviate this difficulty, we propose the use of an implicitly-defined penalty parameter via the DNN approximation of the PDE solution at the previous time-step. We justify our choice of penalty by performing a coercivity analysis of the respective loss functional. Also, for the JKO-type gradient flow, we consider a Sinkhorn-Knopp approximation of the $2$-Wasserstein distance and we find that the resulting ResNET for every time-step is approximating the exact solution to an increasing accuracy concerning the spatial dimension $d$.

The remainder of this work is structured as follows. In Section \ref{pde} we discuss the basic model problem and the (continuous) variational interpretations we shall be concerned with for the design of the approximation algorithms. In Section \ref{numer_ANN}, we discuss various approaches for the design of DNN-based numerical methods and we present the Nitsche functional idea, which is, in turn, used in Section \ref{meat} for the design of the \emph{discrete gradient flow deep Nitsche method} proposed in this work. A series of numerical experiments are presented in Section \ref{num_ex}.  Finally, in Section \ref{JKO_sec}, we present a discrete JKO-type gradient flow for the heat problem, which applies also to other Fokker-Planck type equations, before drawing some conclusions in Section \ref{conclusions}.


\section{Model problem and variational interpretations}\label{pde}

In the functional analytic setting of a Gel'fand triple $\mathcal{V}\subset \mathcal{H}\subset\mathcal{V}^*$, we are primarily concerned with finding an approximate solution $u: (0,T]\times \Omega\to \mathbb{R}$, for a domain $\Omega \subset \mathbb{R}^d $, $d\in\mathbb{N}$, and a $T>0$, satisfying
\begin{equation} \label{equation:1}
	\begin{aligned}
		\partial_t u + \mathcal{L} u=&\ F, \quad \quad \quad \text{in } (0,T] \times \Omega \\
		u(0,\cdot) =&\ u_0, \,   \qquad \text{  in } \Omega \\
		 u= &\ g_D^{}, \quad \quad \quad   \text{on } (0,T] \times \Gamma_D^{},\\
	\mathbf{n}\cdot \nabla u= &\ g_N^{}, \quad \quad \quad   \text{on } (0,T] \times \Gamma_N^{},
	\end{aligned}
\end{equation}
with $\mathcal{L}:\mathcal{V}\to\mathcal{V}^*$ denoting a second-order elliptic self-adjoint partial differential operator admitting a variational minimization interpretation, $\mathbf{n}$ denoting the unit outward normal vector to $\partial\Omega$, $\Gamma_D\cup\Gamma_N=\partial\Omega$, the Dirichlet and Neumann parts of the boundary, respectively. We assume that $u_0\in\mathcal{H}$, $g\in\partial\mathcal{H}$ is sufficiently smooth to be well-defined in the sense of traces in a suitable space $\partial\mathcal{H}$ in what follows, and $F\in L_2(0, T;\mathcal{V}^*)$, in the usual notation of B\"ochner spaces. We stress, nevertheless, that $\mathcal{L}$ being self-adjoint and/or of second order is not an essential restriction for the developments concerning $L_2$-gradient flows below. For simplicity of the exposition, however, we shall confine ourselves to the `canonical' cases $\mathcal{H}=L_2(\Omega)$, $\partial\mathcal{H}=L_2(\partial\Omega)$, and $\mathcal{V}= H^1(\Omega)$ or $\mathcal{V}=H^1_0(\Omega)$ (using standard notation for Lebesgue and Hilbertian Sobolev spaces) and $\mathcal{L}=-\Delta$, with $\Delta$ denoting the Laplacian with respect to the spatial variables. The latter choices give rise to the classical heat problem with various boundary conditions.

Upon considering a subdivision of the time interval $[0,T]$ into timesteps $I_k:=(t_{k-1},t_k]$, $k=1,\dots,N_t$, for $0=:t_0<t_1<\dots<t_{N_t}:=T$, for $N_t\in\mathbb{N}$, \eqref{equation:1} can be discretized in time as follows: let $u^0:=u_0$ and, for $k=1,\dots, N_t$, seek approximations $u^{k}$ such that
\begin{equation}\label{eq:backward_Euler}
	\frac{ u^k-u^{k-1}}{\tau_k} +\mathcal{L}u^k = F^k,
\end{equation}
with $\tau_k:=t_k-t_{k-1}$, and $F^k:=F(t_k,\cdot)$, together with the boundary conditions to ensure closure.

Now, assume that there exists an energy functional $I:\mathcal{V}\to\mathbb{R}$ whose respective Euler-Lagrange operator is given by $\mathcal{L}$. In the `canonical' case, the classical Dirichlet energy 
\begin{equation}\label{eq:dirichlet}
	I(w)=\frac{1}{2}\int_\Omega |\nabla w|^2 \ud x
\end{equation}
corresponds to the Euler-Lagrange operator $\mathcal{L}=-\Delta$ when $\alpha=1$, $\beta=0$, and $g=0$. Moreover, the minimizer 
\[
w=\arg\min_{w\in \mathcal{V}} \Big(I(w)-\int_\Omega F w\ud x\Big)
\]
within the space $\mathcal{V}=H^1_0(\Omega)$ is unique from the Lax-Milgram Lemma. 

An important, yet now classical, observation for what follows is that \eqref{eq:backward_Euler} is the Euler-Lagrange equation of the following variational minimization problem as follows: for $u^0:=u_0$ and $k=1,\dots, N_t$, seek approximations $u^{k}$ defined by
\begin{equation}\label{eq:L2_grad}
u^k=\arg\min_{w\in \mathcal{V}} \bigg(\frac{1}{2}\|w-u^{k-1}\|_{L_2(\Omega)}^2+ \tau_k \Big(I(w)-\int_\Omega  F^k w\ud x\Big)\bigg).
\end{equation}
Thus, instead of solving \eqref{eq:backward_Euler} directly, we can instead construct numerical methods based on solving \eqref{eq:L2_grad} instead. 

Another possibility to write a variational minimization problem for \eqref{eq:backward_Euler} for the particular case $\alpha=1$, $\beta=0$, $g=0$ and $\Omega=\mathbb{R}^d$, is the celebrated JKO-functional introduced by Jordan, Kinderlehrer and Otto in \cite{jko} for the case where the solution $u$ (and the initial condition $u_0$) refer to probability density functions: find  
 \begin{equation}\label{eq:jko}
u^k=\arg\min_{p\in K} \bigg(	\frac{1}{2}\mathcal{W}(u^{k-1},p)^2 + \tau_k \int_{\mathbb{R}^d} p\log p \ud x\bigg),
\end{equation}
over all $p \in K$, where $K$ is the set of all probability densities on $\mathbb{R}^d$ having finite second moments; here $\mathcal{W}$ represents the $2$-\emph{Wasserstein metric} (or \emph{distance}) $\mathcal{W}(\mu_1,\mu_2)$ between two probability measures $\mu_1,\mu_2$ on $\mathbb{R}^d$ given by
\begin{equation}
\mathcal{W}(\mu_1,\mu_2) =\bigg( \inf \limits_{p\in P(\mu_1,\mu_2)} \int_{\mathbb{R}^d \times \mathbb{R}^d } |x-y|^2 p(\!\ud x,\ud y)\bigg)^{1/2};
\end{equation} 
here $ P(\mu_1,\mu_2)$ is the set of probability measures on $\mathbb{R}^d \times \mathbb{R}^d$ with marginals $\mu_1,\mu_2$. (A probability measure $p$ is in $P(\mu_1,\mu_2)$ if and only if for each Borel subset $A \subset \mathbb{R}^d$, $p(A\times \mathbb{R}^d ) = \mu_1(A)$ and $p( \mathbb{R}^d\times A) = \mu_2(A)$ hold.)  It is proved in \cite{jko} that, given $u^0 \in K$, there exists a unique solution of the scheme \eqref{eq:jko}. Moreover, it is proved that the energy functional is convex and that an appropriate interpolation of the solution to \eqref{eq:jko}  for $k=1,\dots,N_t$ converges to the unique solution of the heat equation. We note that the result in \cite{jko} holds for a more general class of Fokker-Planck type equations by suitably modifying the entropy. 


An alternative way of defining the Wasserstein distance is as the infimum taken over all random variables $X, Y$ such that $X$ has distribution $\mu_1$ and $Y$ has distribution $\mu_2$, which is the infimum over all possible couplings of the random variables $X$ and $Y$, namely,
\begin{align}
	\mathcal{W}^2(\mu_1,\mu_2) = \inf \mathbb{E}[|X-Y|^2].
\end{align}

The above variational principles are possible upon understanding subtle properties of the diffusion operators involved, in conjunction with respective boundary conditions. Alternatively, one can always consider the \emph{space-time} `least squares' minimization problem related to \eqref{equation:1}:
\begin{equation}\label{LS}
u=\arg\min_{w\in L_2(0,T;\mathcal{V})}\bigg( || \partial_t w + \mathcal{L}w -F||^2_{L_2(0,T;\mathcal{H})} + || \alpha w +\beta\mathbf{n}\cdot \nabla w- g ||^2_{L_2(0,T;\partial \mathcal{H})} + || w(0,\cdot) - u_0 ||^2_{\mathcal{H}}\bigg),
\end{equation}
which can be used for general classes of PDE problems. We remark that $\mathcal{H}$ here can accept weighted inner products; see \cite{sirignano_dgm_2018} for details.

Although we focus on the case of the heat problem for simplicity of the exposition, we stress that more general cases of PDE problems in the form of  \eqref{equation:1} can be treated by the methodologies presented below.  For instance, any dissipative PDE whose spatial operator is the Euler-Lagrange operator of a suitable energy functional $I(\cdot)$ is admissible for representation by \eqref{eq:L2_grad}, while the JKO theory applies to a general family of linear Fokker-Planck equations also, cf., \cite{jko}.  

\section{Numerical solution of PDE problems via ANNs}\label{numer_ANN}
Perhaps the most common paradigm in employing artificial neural networks  (ANNs) for the numerical solution of PDE problems is the collocation one, whereby the solution sought is modeled by an appropriate ANN architecture and the ANN parameters are optimized to approximate the exact solution through suitably defined loss functionals. 

The recent, yet already highly popular, methodologies of the \emph{Deep Galerkin Method} (DGM)  of Sirignano and Spiliopoulos in \cite{sirignano_dgm_2018} and the original \emph{Physics Informed Neural Network} (PINN) approach of Raissi, Perdikaris \& Karniadakis \cite{PINN} employ various network architectures to model the approximate solution and minimize through the `least-squares' space-time variational problem \eqref{LS}. These methods are very general in terms of the problems to which they can be applied. At the same time, DGM and PINN approaches may not always capture the fine properties of the respective PDE problems, as they are confined to performing minimization within certain function spaces only. In contrast, using more specialized energy functionals, such as the Dirichlet energy \eqref{eq:dirichlet} one restricts the class of problems but allows for minimization in more relevant solution spaces for the problem at hand.

For example, if we consider the elementary case of an (elliptic) Poisson problem (i.e., no time variable) with essential boundary conditions:
\begin{equation}\label{eq:Poisson}
-\Delta u = F \quad\text{in}\quad \Omega ,\qquad u=g\quad\text{on}\quad \partial\Omega
\end{equation} 
within DGM or the PINN approach of \cite{PINN}, it is proposed to find an ANN $U$ (of some given architecture) by minimizing the loss functional
\[
|| \Delta U +F||^2_{L_2(\Omega)} + ||  U - g ||^2_{L_2(\partial\Omega)}\to \min
\]
i.e., minimization of the PDE residual takes place in the $L_2$-norm and the boundary conditions are also enforced in the corresponding sense.  In contrast, using the Dirichlet energy instead, we can compute another ANN approximation $\tilde{U}$ by minimizing:
\begin{equation}\label{eq:dirichlet_energy_pde}
I(\tilde{U})-\int_\Omega F \tilde{U}\ud x =\int_\Omega \Big(\frac{1}{2} |\nabla \tilde{U}|^2- F \tilde{U}\Big)\ud x \to \min
\end{equation}
for $\tilde{U}=g$ on $\partial \Omega$. Standard variational calculus arguments show that $\tilde{U}$ is aiming to minimize the
error $|| \nabla (u-\tilde{U})||_{L_2(\Omega)}$ or, equivalently, the
 functional 
\[
|| \nabla (u-\tilde{U}) ||_{L_2(\Omega)}^2=|| \Delta \tilde{U} +F||^2_{H^{-1}(\Omega)} \to \min
\]
under the condition $\tilde{U}=g$ on $\partial \Omega$, i.e., minimization takes place in a weaker norm. This may be advantageous in many situations of interest and/or higher-dimensional problems. The strategy of considering optimization of ANNs under \eqref{eq:dirichlet_energy_pde} constitutes the, so-called, \emph{Deep Ritz Method} of E and Yu \cite{e_deep_2017}. A key practical challenge for the deep Ritz method is the requirement for the strong imposition of essential boundary conditions. Strongly imposing Dirichlet-type conditions appear to be introducing some stiffness in simulations due to the global nature of the ANN approximants.

The issue of imposition of essential boundary conditions in energy minimization problems of the type \eqref{eq:dirichlet_energy_pde} can be addressed by altering the loss functional to facilitate the weak imposition of essential boundary conditions. More specifically, by seeking ANN minimizer $\hat{U}$ of the modified Dirichlet energy functional
\begin{equation}\label{eq:Nitsche_energy}
\int_\Omega \Big(\frac{1}{2} |\nabla \hat{U}|^2- F \hat{U}\Big)\ud x  -\int_{\partial\Omega} \mathbf{n}\cdot\nabla \hat{U} (\hat{U}-g) \ud s+\int_{\partial\Omega} \frac{\gamma}{2}(\hat{U}-g)^2 \ud s\to \min
\end{equation}
for $\hat{U}\in H^1(\Omega)$, involving a so-called, \emph{penalty parameter} $\gamma$, which should be defined by the user. Notice that when $\hat{U}=g$ the last two terms of the last functional disappear. The second term on the right-hand side of the functional in \eqref{eq:Nitsche_energy} is \emph{not} sign-definite; yet its inclusion is important to retain the consistency of \eqref{eq:Nitsche_energy} with respect to the PDE problem \eqref{eq:Poisson}. To facilitate convexity of the functional, at least in finite-dimensional subspaces of $H^1(\Omega)$, the third term in \eqref{eq:Poisson} is included. To ensure successful `stabilization' (i.e., coercivity in a suitable norm), the penalty parameter $\gamma$ has to be chosen appropriately; this means in practical terms that it should be chosen large enough to `stabilize' effectively but not too large to introduce stiffness in the numerical method. Note that, formally, we expect $\hat{U}\to \tilde{U}$ as $\gamma\to\infty$.

\section{A discrete gradient flow deep Nitsche method}\label{meat}

Approximating an elliptic boundary value problem by minimizing \eqref{eq:Nitsche_energy} for a residual-type network $\hat{U}$ was discussed in \cite{liao_deep_2019}; this method is referred to as the \emph{deep Nitsche method}. To carry out the error analysis in \cite{liao_deep_2019}, the validity of a Bernstein-type inequality (trace inverse estimate) for ANNs was assumed which, in turn, determines the penalty parameter $\gamma$. Given the parameter-dependent growth of ANNs, it is not clear if such an assumption is realistic or not.  Here, we will propose a practical rectification of this state of affairs, which can give a concrete choice of the penalty parameter $\gamma$, which appears to be effective in stabilizing the method in all numerical examples we have performed.

Further, the new practical version of minimization using \eqref{eq:Nitsche_energy} mentioned above will be combined with \eqref{eq:L2_grad} to establish a new ANN  framework, based on deep residual-type architectures, for the approximation of solutions to \eqref{equation:1}. The approximate solution will be, therefore, expressed as a sequence of ANNs $\hat{U}^k$, each representing the approximate solution in the time interval $I_k$. This point of view has several potential salient features compared to the generic `monolithic' space-time residual approximation proposed in \cite{sirignano_dgm_2018, PINN}. In particular, the minimization takes place in the natural $L_2(H^1)$-norm setting, and the boundary conditions are treated weakly to alleviate numerical stiffness. Furthermore, starting from fitting an ANN approximation of the known initial condition of the PDE problem, not involving differentiation of the network, the method progresses iteratively in computing  $\hat{U}^k$, $k=1,\dots, N_t$, having the same architecture as the initial condition ANN. As such, we can use the ANN parameters fitted in the previous timestep as an initial guess for the optimization in the current one; this leads to a significant reduction in epochs needed to achieve a given accuracy. The ANN architecture used in the numerical experiments follows the residual architecture proposed by Sirignano and Spiliopoulos in \cite{sirignano_dgm_2018}.

\subsection{Method definition}

We consider the PDE problem: find $u\in L_2(0,T;H^1(\Omega))$ such that:
\begin{equation} \label{equation:2}
	\begin{aligned}
		\partial_t u -\nabla\cdot (A\nabla u)=&\ F, &\quad \quad \quad \text{in } (0,T] \times \Omega, \\
		u(0,\cdot) =&\ u_0, \,   &\quad \quad \quad \text{   in } \Omega, \\
	 u= &\ g_D^{}, & \quad \quad \quad  \text{on } (0,T] \times \Gamma_D^{},\\
	\mathbf{n}\cdot A\nabla u= &\ g_N^{}, & \quad \quad \quad   \text{on } (0,T] \times \Gamma_N^{},
	\end{aligned}
\end{equation}
in the distributional sense, for a symmetric, uniformly positive definite and bounded diffusion tensor $A\in [L^{\infty}(\Omega)]^{d\times d}$, viz., there exist $\lambda_1,\lambda_d>0$ (the smallest and the largest eigenvalue, respectively,) such that 
\begin{equation}\label{eq:uniform}
\lambda_1 |\xi|^2\le \xi^T A\xi \le \lambda_d |\xi|^2,\quad  \text{ for all }\ \xi\in \mathbb{R}^d.
\end{equation}
We select a subdivision of the time-interval $[0,T]$, reading $0=t_0<t_1<\dots<T_{N_t}$ and we seek approximations $u^k$ to $u(t_k,\cdot)$ given by
\begin{equation}\label{eq:L2_grad_concrete}
	u^k=\arg\min_{w\in H^1(\Omega)} \Big(\frac{1}{2}\|w-u^{k-1}\|_{L_2(\Omega)}^2+ \tau_k N(w)\Big),
\end{equation}
with 
\[
N(w):=
	\int_\Omega \Big(\frac{1}{2} |\sqrt{A}\nabla w|^2- F w\Big)\ud x  -\int_{\Gamma_D} \mathbf{n}\cdot A\nabla w(w-g_D^{}) \ud s+\int_{\Gamma_D} \frac{\sigma}{2}(w-g_{D}^{})^2 \ud s-\int_{\Gamma_N}g_{N}^{}w\ud s.
\]
We note carefully, that the second and third terms in the definition of $N$ vanish when $N$ is evaluated at the exact solution of \eqref{equation:2}; their presence, however, becomes important when $w$ is not satisfying the boundary condition exactly. In particular, the second term ensures the consistency of the functional with respect to the  PDE problem \eqref{equation:2}, while the third term penalizes $w$ when its value on the Dirichlet part of the boundary $\Gamma_D$ differs from the Dirichlet boundary condition. A similar in spirit, functional has appeared within the context of minimization problems in the presence of interfaces in \cite{makridakis}. We stress that more general classes of problems can be treated in this context, e.g., nonlinear reaction-diffusion equations with mixed boundary conditions, etc. We refrain from providing the full generality here to focus on the key mathematical and algorithmic issues.

A crucial aspect of any deep learning based algorithm is its \emph{architecture}: clever choices of architectures, may help exploit a prior knowledge about the application and therefore improve the performance.  In \cite{sirignano_dgm_2018}, the authors propose an interesting architecture, similar to Recurrent Neural Networks (RNNs) and, more specifically, similar to the Long Short-Term Memory (LSTM) architecture \cite{hochreiter_long_1997}. LSTM networks are well-suited for time (and space) dependent data, since there can be lags of unknown duration between important events. Furthermore, LSTMs deal with the vanishing gradient problem that can be encountered when training traditional RNNs \cite{DBLP:journals/corr/abs-1211-5063}. The shape of the solution $u(t, x)$ for $t<T$, although smooth, to the presence of diffusion, may rapidly change in certain spatial regions. For $d_{in}$ as input dimension and $d_{out}$ as output dimension, we consider  $\bm{x} \in \mathbb{R}^{d_{in}}$, $W^1 \in \mathbb{R}^{n_1\times d_{in}}, b^1 \in \mathbb{R}^{n_1}$, $U^{\cdot,l} \in \mathbb{R}^{n_{l+1}\times d_{in}}$, $W^{\cdot,l} \in \mathbb{R}^{n_{l+1}\times n_l}$,$b^{\cdot,l} \in \mathbb{R}^{n_{l}}$, for $l=1,...,L$, and $W \in \mathbb{R}^{d_{out},n_{L+1}}, b \in \mathbb{R}^{d_{out}}$, where $n_l$ is the dimension of $S_l$ that is equal to the number of neurons at each block. The architecture proposed in \cite{sirignano_dgm_2018} reads:
\begin{equation} \label{eq:resnet_one}
	\begin{aligned}
		S^1 &=\sigma(W^1\bm{x}+b^1) \\
		Z^l &= \sigma(U^{z,l}\bm{x} + W^{z,l}S^l +b^{z,l}), \quad l=1,..,L \\
		G^l &= \sigma(U^{g,l}\bm{x} + W^{g,l}S^l +b^{g,l}), \quad l=1,..,L \\
		R^l &= \sigma(U^{r,l}\bm{x} + W^{r,l}S^l +b^{r,l}), \quad l=1,..,L \\
		H^l &= \sigma(U^{h,l}\bm{x} + W^{h,l}(S^l \odot R^l) +b^{h,l}), \quad l=1,..,L	\\
		S^{l+1} &= (1-G^l)\odot H^l + Z^l \odot S^l, \quad l=1,..,L \\
		f(t,\bm{x};\theta)&=WS^{L+1}+b 
	\end{aligned}
\end{equation}

\noindent where $\bm{x} = (t,x)$, the number of blocks are $L+1$ and $\odot$ is the element-wise multiplication operator. The parameters of this architecture are:
\[
\theta = \big\{W^1,b^1, \big(U^{z,l},W^{z,l},b^{z,l}  \big)^L_{l=1}, \big(U^{g,l},W^{g,l},b^{g,l}  \big)^L_{l=1}, \big(U^{r,l},W^{r,l},b^{r,l}  \big)^L_{l=1}, \big(U^{h,l},W^{h,l},b^{h,l}  \big)^L_{l=1},W,b  \big\}
\]

Each block parameter consists of $M$ neurons/units with element-wise non-linearity function $\sigma: \mathbb{R}^M \to \mathbb{R}^M$ defined as $\sigma(z) = \big(\phi(z_1),...,\phi(z_M)\big)$ where $\phi: \mathbb{R} \to \mathbb{R}$ is the activation function. 

The aforementioned mathematical description of the architecture appears to be complicated at first sight. Each layer takes as an input the batch inputs $\bm{x}$ (for instance, a set of randomly sampled time-space points) and the output $S^l$ of the previous layer/block. The final block $S^{L+1}$ passes through a linear transformation to produce the final output $y=f(t,\bm{x};\theta)$. Compared to a Multilayer Perceptron (MLP) of the same number of neurons/units, the number of parameters in each hidden layer of the above architecture is approximately eight times greater than the same number in a common dense layer. This observation derives from the fact that each layer in \eqref{eq:resnet_one} has $8$ weight matrices and $4$ bias vectors.

\subsection{The choice of the penalty parameter $\gamma$}
A crucial aspect of the success of the proposed method is the judicious choice of the penalty parameter $\gamma$. In the context of finite element methods with piecewise polynomial approximants on triangulations and Nitsche-type imposition of boundary conditions, the choice of the penalty parameter is determined by the constant of, so-called, \emph{inverse estimates}. To showcase the use of such estimates, we consider momentarily a triangulation $\mathcal{T}(\mathcal{N})$ of $\Omega$ into mutually non-overlapping $d$-dimensional simplices with $\mathcal{N}$ denoting the set of nodes (vertices) with respect to which the triangulation is defined. (If $\Omega$ has curved boundaries we approximate the domain, as is standard in finite element methods.) Let also $S_h$ be the space of continuous piecewise polynomial functions of degree at most $p$, subordinate to the triangulation.  Then, there exists a constant $C>0$, independent of the triangulation granularity, (depending though on the polynomial degree and, possibly on elemental shapes,) such that the estimate 
\begin{equation}\label{eq:inv_in_poly}
\|\sqrt{h}v\|_{L^2(\Gamma_D)} \leq C\|v\|_{L^2(\Omega)},
\end{equation}
holds, with $h$ the local mesh function, representing the local maximum distance between nodes in the triangulation. In particular, in the vicinity of the boundary $\Gamma_D$, $h$ is proportional to the distance of the locally nearest non-boundary node from  $\Gamma_D$.  We now describe why inverse estimates such as \eqref{eq:inv_in_poly} prescribe suitable choices for $\sigma$, if we seek to minimize \eqref{eq:L2_grad_concrete} within $S_h$. To that end, we have, respectively, for $u_h\in S_h$,
\[
\begin{aligned}
N(u_h)
\ge&\ 
\int_\Omega \Big(\frac{1}{2} |\sqrt{A}\nabla u_h|^2- F u_h\Big)\ud x  -
\|\sqrt{A h}\nabla u_h\|_{L^2(\Gamma_D)}
\||\sqrt{A}\mathbf{n}|h^{-1/2}(u_h-g_D^{})\|_{L^2(\Gamma_D)} \\
&\quad +\int_{\Gamma_D} \frac{\sigma}{2}(u_h-g_{D}^{})^2 \ud s-\int_{\Gamma_N}g_{N}^{}u_h\ud s\\
\ge&\ 
\int_\Omega \Big(\frac{1}{4} |\sqrt{A}\nabla u_h|^2- F u_h\Big)\ud x  
 +\frac{1}{2}\int_{\Gamma_D} \Big(\sigma-\frac{2C^2\lambda_d^2}{\lambda_1 h}\Big)(u_h-g_{D}^{})^2 \ud s-\int_{\Gamma_N}g_{N}^{}u_h\ud s,
\end{aligned}
\]
using \eqref{eq:uniform} and \eqref{eq:inv_in_poly} by observing that $\nabla u_h$ is still a piecewise polynomial over the same triangulation, along with the inequality $\alpha\beta\le \alpha^2/4+\beta^2$ for $\alpha,\beta\in\mathbb{R}$ in the last step. Therefore, to ensure the positivity of the quadratic terms and, therefore, uniqueness of the solution, we should select $\sigma\ge 2C^2\lambda_d^2/(\lambda_1 h)$.  At the other end of the spectrum, the penalty parameter $\sigma$ should \emph{not} be chosen `too' large, as this will effectively annihilate the advantages brought by the use of weak imposition of boundary conditions.

Unfortunately, inverse estimates such as \eqref{eq:inv_in_poly} are not available for general (deep) ANNs and, therefore, it is not \emph{a priori} clear how to select $\gamma$ in \eqref{eq:L2_grad_concrete} for ANN functions. Nevertheless, a good choice of the penalty parameter $\gamma$ is essential for the success of the proposed method. To address this we now provide a `semi-heuristic' construction to facilitate a practical choice for $\gamma$.  Following the standard practice, in the implementation below, we approximate the integrals in the minimization principles by quadrature rules based on random point clouds $\mathcal{N}_I:=\{x_n^I\}_{n=1}^{\mathrm{N}_I}$ contained in $ \Omega$ according to a uniform distribution, $\mathcal{N}_D:=\{x_n^D\}_{n=1}^{\mathrm{N}_D}\subset \Gamma_D$ and  $\mathcal{N}_N:=\{x_n^N\}_{n=1}^{\mathrm{N}_N}$; in practice, we select the uniform distribution, unless there is a specific modeling reason to make a different choice.  We note that one of the two boundary point clouds may be void if the respective boundary condition is not enforced. For uniformly distributed point cloud distributions, we approximate the integrals, as is standard, in a Monte Carlo fashion:
\[
\int_\Omega f \ud x \approx \frac{|\Omega|}{\mathrm{N}_I}\sum_{n=1}^{\mathrm{N}_I} f(x_n^I),\qquad \int_{\Gamma_D} g \ud s \approx \frac{|\Gamma_D|}{\mathrm{N}_D}\sum_{n=1}^{\mathrm{N}_D} g(x_n^D),
\qquad \int_{\Gamma_N} g  \ud s \approx \frac{|\Gamma_N|}{\mathrm{N}_N}\sum_{n=1}^{\mathrm{N}_N} g(x_n^N).
\]

Using the above, we define the approximate Nitsche functional $\tilde{N}(w,w)\approx N(w)$, given by
\begin{equation}
\begin{aligned}
\tilde{N}(w;v):=&\ 
\frac{|\Omega|}{\mathrm{N}_I}\sum_{n=1}^{\mathrm{N}_I} \Big(\frac{1}{2} |\sqrt{A}(x_n^I)\nabla w(x_n^I)|^2- F(x_n^I) w(x_n^I)\Big) \\
& -\frac{|\Gamma_D|}{\mathrm{N}_D}\sum_{n=1}^{\mathrm{N}_D} \big(\mathbf{n}(x_n^D)\cdot A(x_n^D)\nabla w(x_n^D)(w(x_n^D)-g_D^{}(x_n^D))\big)\\
&+\frac{|\Gamma_D|}{\mathrm{N}_D}\sum_{n=1}^{\mathrm{N}_D} \frac{\gamma(x_n^D;v)}{2}(w(x_n^D)-g_{D}^{}(x_n^D))^2 -\frac{|\Gamma_N|}{\mathrm{N}_N}\sum_{n=1}^{\mathrm{N}_N} g_{N}^{}(x_n^N)w(x_n^N),
\end{aligned}
\end{equation}
and we seek to minimize instead the loss functional given by \eqref{eq:L2_grad_concrete} with $\tilde{N}$ replacing $N$. Note that we allow the penalty $\gamma$ to depend independently on an ANN $v$, which may differ from $w$; this is non-standard and will be discussed in detail below.

A key attribute of the proposed approach is the estimation of the penalty parameter $\sigma$, which we now consider. Standard estimation of the indefinite term yields
\begin{equation}\label{inlieuof}
\begin{aligned}
	& \Big|\frac{|\Gamma_D|}{\mathrm{N}_D}\sum_{n=1}^{\mathrm{N}_D} \big(\mathbf{n}(x_n^D)\cdot A(x_n^D)\nabla w(x_n^D)(w(x_n^D)-g_D^{}(x_n^D))\big) \Big|\\
		\le &\ \frac{|\Omega|}{\mathrm{N}_I}\sum_{n=1}^{\mathrm{N}_D}  \frac{|\Gamma_D|\mathrm{N}_I \lambda_d^2}{|\Omega|\mathrm{N}_D\tau(x_n^D)}|\nabla w(x_n^D)|^2 + \frac{|\Gamma_D|}{\mathrm{N}_D}\sum_{n=1}^{\mathrm{N}_D}\frac{\tau(x_n^D)}{4}(w(x_n^D)-g_D^{}(x_n^D))^2,
\end{aligned}
\end{equation}
for any $\tau(x_n^D)>0$, $x_n^D\in \mathcal{N}_D$. Now, let $y_n\equiv y_n(x_n^D)\in \mathcal{N}_I$ denote the nearest point to $x_n^D\in \mathcal{N}_D$, such that $|\nabla w(y_n)|\neq 0$, that has not been selected in a previous boundary point. 
Combining \eqref{inlieuof}, along with the uniform ellipticity, and selecting
\begin{equation}\label{pen}
\gamma(x_n^D;w)\ge \tau(x_n^D):=8\frac{|\Gamma_D|\mathrm{N}_I \lambda_d^2}{|\Omega|\mathrm{N}_D\lambda_1}\frac{|\nabla w(x_n^D)|^2 }{|\nabla w(y_n)|^2 }, 
\end{equation}
noting the dependence of $\gamma$ on $w$, we deduce 
\begin{equation}\label{coer_discrete}
\begin{aligned}
\tilde{N}(w;w)\ge &\ 
	\frac{|\Omega|}{\mathrm{N}_I}\sum_{n=1}^{\mathrm{N}_I} \Big(\frac{1}{4} |\sqrt{A}(x_n^I)\nabla w(x_n^I)|^2- F(x_n^I) w(x_n^I)\Big) \\
	&+\frac{|\Gamma_D|}{\mathrm{N}_D}\sum_{n=1}^{\mathrm{N}_D} \frac{\gamma(x_n^D;w)}{4}(w(x_n^D)-g_{D}^{}(x_n^D))^2 -\frac{|\Gamma_N|}{\mathrm{N}_N}\sum_{n=1}^{\mathrm{N}_N} g_{N}^{}(x_n^N)w(x_n^N), 
\end{aligned}
\end{equation}
with the second argument in $\tilde{N}$ referring to the (nonlinear) dependence of $\gamma$ on $w$.

Therefore, the coercivity  of $\tilde{N}(w;w)$ is possible upon selecting the penalty $\gamma$ according to \eqref{pen}. Note that, as it stands, the choice of $\gamma$ is \emph{nonlinear} with respect to the parameters $\theta$. This is non-standard in Nitsche-type/weak imposition of boundary conditions. Here, however, we use this concept to introduce some quantitative information on the size of $\gamma$, which is crucial for the success of the method as we shall see below. In a practical method, to avoid back-propagating the additional $w$-nonlinearity introduced in the penalty $\gamma$, in the gradient descent step we use the value $u^k(\theta_{m-1})$ to evaluate the penalty instead, viz., we implement $\tilde{N}(u^k(\theta_m);u^k(\theta_{m-1}))$ at the current step; we refer to Algorithm \ref{alg:time-dnml} for details of the implemented process. Finally, we note that the values $\theta_0$ used are the ones from the previous time-step.

\section{Numerical experiments}\label{num_ex}

We will now present a comprehensive numerical experiment showcasing the performance of the proposed ``Nitsche-type'' method used to approximate solutions to linear heat-type parabolic problems for $d\in[2,20]$. More specifically, we will demonstrate the good performance of the combined use of the functional $\tilde{N}$ above with the choice of penalty given by \eqref{pen} in conjunction with the architecture described in \eqref{eq:resnet_one}.

We evaluate the performance of the time deep Nitsche method where for $k=1,\dots, N_t$ we compute ANN $u^k$ using
\begin{equation}\label{eq:L2_grad_concrete_discrete}
	u^k=\arg\min_{w\in H^1(\Omega)} \bigg(\frac{|\Omega|}{2\mathrm{N}_I}\sum_{n=1}^{\mathrm{N}_I} (w-u^{k-1})^2(x_n^I)+ \tau_k \tilde{N}(w;w)\bigg).
\end{equation}
Therefore, in contrast to the DGM/PINN approach \cite{sirignano_dgm_2018, jentzen, PINN}, whereby a unique `space-time' ANN/ResNET is computed, the present method constructs one ResNET at each time-step. As a result, the computed ``time-instance" ResNETs can be chosen to have a significantly smaller width and depth, compared to full space-time collocation approaches, resulting in faster performance. Moreover, the optimized ResNET parameters at time-step $k-1$ can be used as excellent initializations for the next timestep $k$, thereby reducing significantly the number of epochs required. Note that this process starts at $k=0$ with the approximation of the initial condition, i.e., a simple function fitting step not involving differentiation of the ResNET; this is typically achieved to very high accuracy. As a result, we expect the approximation $u^1$ to be also good, and so on.  The pseudocode of the algorithm for the time deep Nitsche method is given in Algorithm \ref{alg:time-dnml}.

\begin{algorithm}
\caption{The Time Deep Nitsche Method}
\begin{algorithmic}[1]
\STATE Initialize randomly the parameter set $\theta_0$ and set the learning rate schedule $a_n$. Initialize time step $\tau$.
\STATE Initialize ANN approximating the initial condition $u^0=\arg\min_{w\in H^1(\Omega)}\|w-u_{0}\|^2_{L_2(\Omega)}$.
\FOR{$k=1,2,\dots,N_t$}
\STATE Create ANN $u^k$ equal to $u^{k-1}$ from the previous time step.
    \FOR{each epoch $m=1,2,\dots$}
    \STATE Generate  point clouds $\mathcal{N}_I:=\{x_n^I\}_{n=1}^{\mathrm{N}_I}\subset \Omega$,  $\mathcal{N}_D:=\{x_n^D\}_{n=1}^{\mathrm{N}_D}\subset \Gamma_D$, and  $\mathcal{N}_N:=\{x_n^N\}_{n=1}^{\mathrm{N}_N}\subset\Gamma_N$.
    \STATE Calculate $L(\theta_m):=\frac{|\Omega|}{\mathrm{N}_I}\sum_{n=1}^{\mathrm{N}_I} (u^k(\theta_m)-u^{k-1})^2(x_n^I)+ \tau_k \tilde{N}(u^k(\theta_m);u^k(\theta_{m-1}))$.
    \STATE Take a descent step at the random batch of generated points using SGD: $\theta_{m+1} = \theta_m - a_n \nabla_{\theta} L(\theta_m) $.
    \STATE Exit when $||\theta_{m+1}-\theta_m||$ is small enough.
    \ENDFOR
\ENDFOR
\STATE Return $u^k$ solution for every time step.
\end{algorithmic}
\label{alg:time-dnml}
\end{algorithm}

			
			
				
					
To test the Deep Nitsche Method, we consider \eqref{equation:2} with $A=I_{d\times d}$, and $u_0$ and $f$ so that
\begin{align*}
	u(t,\bm{x}) = \sin(t) \sin (x_1+x_2+...+x_d) \numberthis
\end{align*}
is the exact solution, with  $\bm{x}=[x_1,x_2,...,x_d]^T \in [0,1]^{d}=:\Omega$, for $t\in(0,1]$ and $d$ denotes the spatial dimensions; we take $d = 2,3,5,10$ and $20$. Finally, we consider a uniform time-step $\tau = 0.01$.

To realize the method, we uniformly sample $600 d$ points from the domain and $600$ from each boundary edge, using $3$ blocks with $50$ neurons each, and with $ReLU(\cdot)$ activation. We select the penalty at each iteration $k$ to be
\begin{align}
    \gamma^k = \max_{x_n^D \in \mathcal{N}_D} 500 \cdot  \frac{|\Gamma_D|\mathrm{N}_I \lambda_d^2}{|\Omega|\mathrm{N}_D\lambda_1}\frac{|\nabla w(x_n^D)|^2 }{|\nabla w(y_n)|^2 }.
\end{align}
The weights of the initial condition are initialized using Xavier initialization \cite{pmlr-v9-glorot10a}. The initial network is trained for $50,000$ epochs for $u^{1}$, using the following learning rate schedule 
\begin{align}
	lr(epochs)=
	\begin{cases}
		10^{-2}, \textit{ 1 $\leq$ epochs $<$ 10,000}\\
		10^{-3}, \textit{ 10,000 $\leq$ epochs $<$ 40,000}\\
		10^{-4},\textit{ epochs $\geq$ 40,000}.
	\end{cases}
\end{align} The subsequent networks $u^k$ are trained via $2,000$ epochs  $k=1,\dots,N_t:=1/\tau$ with a constant learning rate equal to $10^{-3}$ for the first $500$ epochs and $10^{-4}$ for the rest. The performance of the method in the above problem is presented in Table \ref{table:TDNMresults}. We note the `scalability' of the method with respect to the dimension $d$ in terms of errors: the metrics taken for $d=2$ are comparable to metrics taken for $d=20$. 

\begin{table}[ht]
	\centering
	\bgroup
	\def\arraystretch{1.5}
	\begin{tabular}{|c|c|c|c|}
		\hline
		\textbf{$d$} & $L^2$ relative error & max error& mean error \\	\hhline{|=|=|=|=|} 
		2  & $1.6e{-2}$  & $9.3e{-3}$ & $2.9e{-3}$  \\ \hline
		3  & $4.7e{-3}$  & $7.2e{-3}$ & $1.3e{-3}$  \\ \hline
		5  & $2.0 e{-3}$  & $1.5e{-3}$ & $2.9e{-4}$  \\ \hline
		10  & $3.5 e{-3}$  & $1.9e{-3}$ & $2.4e{-4}$  \\ \hline
		20  &  $4.2 e{-3}$  & $3.4e{-3}$ & $3.7e{-4}$  \\ \hline
	\end{tabular}
	\egroup
	\caption{Results of the time deep Nitsche method - Varying Penalty}
	\label{table:TDNMresults}
\end{table}

To assess the potential practical advantages of the proposed time deep Nitsche Method, we compare it against the general purpose DGM approach in which a discrete version of \eqref{LS} is solved for the same PDE problem. For the solution we use uniformly sampled $600 \cdot (d+1)$ points from the domain and $600$ from each boundary edge, using $3$ blocks with $50$ neurons each, along with the $ReLU(\cdot)$ activation function, the weights are initialized via the Xavier initialization \cite{pmlr-v9-glorot10a} and the network is trained for $200,000$ epochs; the above are in accordance to the implementation in \cite{sirignano_dgm_2018}. We use the same learning rate and the same total number of epochs to make a fair comparison. The performance of the DGM approach is given in Table \ref{table:DGM}.

\begin{table}[ht]
	\centering
	\bgroup
	\def\arraystretch{1.5}
	\begin{tabular}{|c|c|c|c|}
		\hline
		\textbf{$d$} & $L^2$ relative error & max error& mean error \\	\hhline{|=|=|=|=|} 
		2  & $9.9e{-2}$  & $8.7e{-2}$ & $3.4 e{-2}$  \\ \hline
		3  & $9.3e{-2}$  & $1.0e{-1}$ & $3.2 e{-2}$  \\ \hline
		5  & $1.2 e{-1}$  & $1.2e{-1}$ & $2.8e{-2}$  \\ \hline
		10  & $4.5 e{-1}$  & $2.9e{-1}$ & $1.5e{-1}$  \\ \hline
		20  &  $4.7 e{-1}$  & $4.2e{-1}$ & $1.3e{-1}$  \\ \hline
	\end{tabular}
	\egroup
	\caption{Results of the Deep Galerkin Method}
	\label{table:DGM}
\end{table} 

By observing the tables we observe that the time deep Nitsche method outperforms DGM by an order of magnitude for this particular PDE problem. We conclude, therefore, that, at least for this particular PDE problem and with the above architectures and method parameters the time deep Nitsche method is at least as competitive as DGM/PINN while using significantly fewer computational resources in the optimization step.

\section{A discrete gradient flow based on the JKO functional}\label{JKO_sec}

Finally, for the case $F=0$, $\Gamma_D=\emptyset$, $g_N^{}=0$, i.e., for the pure homogeneous Neumann problem, we also present results from an implementation of seeking approximations based on minimizing the JKO-functional \eqref{eq:jko}. The challenge of evaluating the $2$-Wasserstein distance is addressed via the Sinkhorn-Knopp Algorithm \cite{sink}. More specifically, we use the algorithm proposed in \cite{sink_algo}. The discretization scheme is the following
 \begin{align}
\frac{1}{2}\mathcal{W}^2(p^{k-1},p)+ h \int_{\Omega} p\log p 
\end{align}

\noindent over the appropriate class of densities $p \in K$ where $K$ is the set of all probability densities in $\Omega$ having finite second moments and $\mathcal{W}(\cdot,\cdot)$ the $2$-Wasserstein distance. Based upon the aforementioned scheme we can derive a learning algorithm for neural networks that solves the heat equation. The algorithm is given below.

\begin{algorithm}
\caption{The Deep Wasserstein Method}
\begin{algorithmic}[1]
\STATE Initialize randomly the parameter set $\theta_0$ and set the learning rate schedule $a_n$. Initialize time step $\tau$.
\STATE Initialize ANN approximating the initial condition $u^0=\arg\min_{w\in H^1(\Omega)}\|w-u_{0}\|^2_{L_2(\Omega)}$.
\FOR{$k=1,2,\dots,N_t$}
\STATE Create ANN $u^k$ equal to $u^{k-1}$ from the previous time step.
    \FOR{each epoch $m=1,2,\dots$}
    \STATE Generate  interior point cloud $\mathcal{N}_I:=\{x_n^I\}_{n=1}^{\mathrm{N}_I}\subset \Omega$.
    \STATE Calculate $L(\theta_m) = \frac{1}{2}d(u^{k-1},u^k(\theta_m))^2 + \tau \int_{\Omega} u^k(\theta_m)\log u^k(\theta_m) $.
    \STATE Take a descent step at the random batch of generated points using SGD: $\theta_{m+1} = \theta_m - a_n \nabla_{\theta} L(\theta_m) $.
    \STATE Exit when $||\theta_{m+1}-\theta_m||$ is small enough.
    \ENDFOR
\ENDFOR
\STATE Return $u^k$ solution for every time step.
\end{algorithmic}
\label{alg:time-wasser}
\end{algorithm}

\par We use our method to solve a parabolic type PDE with Neumann boundary conditions. The problem has the following analytical solution
\begin{align*}
u(t,\bm{x}) = \frac{1}{2}(e^{-d\pi^2t}\prod_{i=1}^d \cos(\pi x_i)+2), \numberthis
\end{align*}

\noindent where $\bm{x}=[x_1,x_2,...,x_d]^T \in [0,1]^{d}$. We solve this problem for $d =2,3,5,10,20,40$ and $50$ dimensions for $t \in [0,T=1]$ and $\tau = h =  0.01$. 
\par We uniformly sample $100 \cdot d$ points from the domain, we used $3$ blocks with $50$ neurons each, we applied as activation function the $ReLU(\cdot)$ and we initialized the weights using Xavier initialization \cite{pmlr-v9-glorot10a}. The initial network is trained for $4,000$ epochs for $u^0$, using the following learning rate schedule 
\begin{align}
	lr(epochs)=
	\begin{cases}
		10^{-3}, \textit{ 1 $\leq$ epochs $<$ 2,000}\\
		10^{-4},\textit{ epochs $\geq$ 2,000}.
	\end{cases}
\end{align} The subsequent networks $u^k$ are trained for $100$ epochs with constant learning rate equal to $10^{-5}$. The performance of this method is given in Table \ref{tab:wasser}. By observing the results, we conclude that the deep Wasserstein method solves the particular PDE efficiently, especially in higher dimensions in relatively few epochs. It would be interesting to also examine the different functionals presented in \cite{new_JKO} using the ResNET architecture employed in this work. 

\begin{table}[ht]
	\centering
	\bgroup
	\def\arraystretch{1.5}
	\begin{tabular}{|c|c|c|c|}
		\hline
		$d$ & $L^2$ Relative Error & Maximum Error & Mean Error \\	\hhline{|=|=|=|=|}
		2  &  8.7e-2  & 1.0e-1 & 8.5e-2   \\ \hline
		3  & 1.7e-1 & 4.8e-1    & 1.3e-1  \\ \hline
		5  & 8.6e-2  & 4.2e-1  & 5.3e-2   \\ \hline
		10  & 8.2e-3  &  3.9e-2 & 6.8e-3  \\ \hline
		20  & 4.0e-3   & 6.8e-3 & 4.0e-3  \\ \hline
		40  &   2.1e-3  & 4.2e-3 & 1.9e-3  \\ \hline
		50  &  2.5e-3  &  1.6e-2 & 2.5e-3   \\ \hline
	\end{tabular}
	\egroup
	\caption{Results of the deep Wasserstein method}
	\label{tab:wasser}
\end{table}


\section{Acknowledgments}
This research work was supported by the Hellenic Foundation for Research and Innovation (H.F.R.I.) under 1) the ``First Call for H.F.R.I. Research Projects to support Faculty members and Researchers and the procurement of high-cost research equipment grant" (Project Numbers: 3270, 1034 and 2152). Also, EHG wishes to acknowledge the financial support of The Leverhulme Trust (grant number RPG-2021-238) and of EPSRC (grant number EP/W005840/1). 

\nocite{*}
\linespread{1.5}
\bibliographystyle{alpha} 


\newcommand{\etalchar}[1]{$^{#1}$}

\end{document}